\documentclass[11pt]{amsart}
\usepackage{amsthm, amssymb, latexsym, amsmath, color}
\usepackage[a4paper, total={6.5in, 8in}]{geometry}
\usepackage[dvipsnames,svgnames,table]{xcolor}
\usepackage[colorlinks=true,linkcolor=RoyalBlue,urlcolor=RoyalBlue,citecolor=PineGreen]{hyperref}
\usepackage{tikz}
\usetikzlibrary{angles,quotes}
\usepackage{tikz-cd}
\usetikzlibrary{arrows.meta}
\usepackage{comment, caption}
\usepackage{graphicx,subcaption}
\usepackage{float}
\usepackage{pgf,tikz,pgfplots}
\pgfplotsset{compat=1.14}
\usepackage{mathrsfs}
\usetikzlibrary{arrows}
\usepackage{capt-of}
\usetikzlibrary{decorations.pathreplacing,cd,positioning,calc}
\usepackage{algorithmic}
\usepgfplotslibrary{fillbetween}
\usetikzlibrary{intersections,patterns}

\usepackage[nameinlink]{cleveref}

\title[Small cells in great-circle arrangements]{Small cells in arrangements of great circles}
\author{Krishnendu Bhowmick}
\address{Krishnendu Bhowmick \\ TU Chemnitz, Germany.}
\email{krish.combi@gmail.com}
\date{}
\hypersetup{
  pdftitle={Small cells in arrangements of great circles},
  pdfauthor={P\'eter \'Agoston and Krishnendu Bhowmick},
  pdfsubject={Uniform bounds for the minimum area and area ratio of great-circle cells}
}

\newcommand{\R}{\mathbb{R}}

\newcommand{\cA}{\mathcal{A}}

\newcommand{\cC}{\mathcal{C}}

\newtheorem{lemma}{Lemma}

\newtheorem{theorem}{Theorem}

\newtheorem{conjecture}{Conjecture}\crefname{conjecture}{conjecture}{conjectures}\Crefname{conjecture}{Conjecture}{Conjectures}

\theoremstyle{remark}
\newtheorem{remark}{Remark}

\begin{document}
\begin{abstract}
We study the areas of the cells determined by an arrangement $\cA$ of
$n$ distinct great circles on the unit sphere, with no three concurrent.
Let $m(\cA)$ and $M(\cA)$ denote the smallest and largest cell areas,
respectively. Our first result is the bound
\[
 m(\cA)\le \frac{C_*+O(n^{-3/4})}{n^2},
 \qquad C_*=7.0666071493\ldots ,
\]
where the implied constant is absolute. This improves the coefficient
$\pi^2$ obtained from classical polygonal perimeter estimates. We also
prove that $M(\cA)/m(\cA)>13$ for all sufficiently large $n$, improving
the bound $7.43$ of L.~Fejes T\'oth. The main tool is a perimeter
inequality for planar polygons which takes account of their exterior
angles. We combine this with a positive semidefinite kernel on the
sphere and use gnomonic projection to obtain bounds for spherical cells.
\end{abstract}
\maketitle
\section{Introduction}

Let $\cA$ be an arrangement of $n\ge3$ distinct great circles on the unit
sphere $S^2\subset\R^3$, with no three circles passing through the same
point. Its cells are convex spherical polygons. Let $\cC(\cA)$ denote
the set of cells, and write
\[
 m(\cA)=\min_{C\in\cC(\cA)}\operatorname{area}(C),\qquad
 M(\cA)=\max_{C\in\cC(\cA)}\operatorname{area}(C),
\]
and put
\[
 H(n)=\sup_{\cA}m(\cA),\qquad
 R(n)=\inf_{\cA}\frac{M(\cA)}{m(\cA)}.
\]
Thus $H(n)$ measures how large the smallest cell can be, while $R(n)$
measures how close the largest and smallest cell areas can be. For
context, the number of cells is $F=n^2-n+2$, and their areas sum to
$4\pi$. It follows immediately that $H(n)\le4\pi/F$.

In 1987, L.~Fejes T\'oth~\cite{FT87} proposed the following conjecture
about the ratio of the largest and smallest cell areas.
\begin{conjecture}[L.~Fejes T\'oth~\cite{FT87}]
\label{conj:ratio}
One has $R(n)\to\infty$ as $n\to\infty$.
\end{conjecture}
In other words, for every fixed $K>0$, every arrangement of sufficiently
many great circles in general position has a cell whose area is more
than $K$ times the area of a smallest cell.

Fejes T\'oth proved that $R(n)>7.43$ for all sufficiently large $n$.
More precisely,
his argument gives, for each fixed $q<Q_0$, the bound $R(n)>q$ for all
sufficiently large $n$, where
\begin{equation}\label{eq:old-ratio}
 Q_0=\left(\frac{8-\pi+4\sqrt{4-\pi}}{\pi}\right)^2
     =7.431834008\ldots .
\end{equation}
The estimate $7.43$ is also recorded in the survey of G.~Fejes T\'oth,
L.~Fejes T\'oth and Kuperberg~\cite[\S19]{FTFK23}.

A stronger small-cell conjecture appears as Conjecture~7 in
Section~11.2 of Brass, Moser and Pach~\cite[p.~446]{BMP05}.
\begin{conjecture}[L.~Fejes T\'oth, as stated in~\cite{BMP05}]
\label{conj:small-cell}
One has $H(n)=o(n^{-2})$ as $n\to\infty$.
\end{conjecture}
Note that \cref{conj:small-cell} implies \cref{conj:ratio}. Indeed,
$M(\cA)\ge4\pi/F$ and $m(\cA)\le H(n)$, so
\[
 R(n)\ge\frac{4\pi}{F H(n)}.
\]
The right side tends to infinity if $H(n)=o(n^{-2})$.
The converse does not follow from this argument, since the largest cell
need not have area comparable to the average.

In this paper, we prove improved bounds for both $H(n)$ and $R(n)$.
We first state the bound for $H(n)$. To describe the constant which
appears in this result, define
\begin{equation}\label{eq:constants}
 \begin{gathered}
 a_* = \sqrt{\pi^2/4-1},\qquad
 d_* = a_*-\arcsin(2a_*/\pi),\qquad
 b=2\cos(d_*/2),\\
 I=\int_0^1(1-x^2)^{-1/4}\,dx
   =\frac12 B\!\left(\frac12,\frac34\right),\qquad
 C_* = \left(\frac{2\pi}{bI}\right)^2.
 \end{gathered}
\end{equation}
Here $B$ is the beta function. Numerically,
$b=1.9727258781\ldots$, $I=1.1981402347\ldots$, and
$C_*=7.0666071493\ldots$.

\begin{theorem}\label{thm:area}
There is an absolute constant $K$ such that, for all sufficiently large
$n$, every arrangement $\cA$ of $n$ distinct great circles on $S^2$, with
no three concurrent, satisfies
\[
 m(\cA)\le\frac{C_*+Kn^{-3/4}}{n^2}.
\]
In particular, $H(n)\le(C_*+O(n^{-3/4}))/n^2$, with $C_*<7.067$.
\end{theorem}

For comparison, the polygonal perimeter estimate used by
Fejes T\'oth already implies
\begin{equation}\label{eq:old-area}
 H(n)\le\frac{\pi^2+o(1)}{n^2}.
\end{equation}
We give the short deduction of~\eqref{eq:old-area} in
\cref{rem:old-area}. This is a consequence of the method in~\cite{FT87},
rather than a bound explicitly stated there. Thus \cref{thm:area}
improves the coefficient from $\pi^2=9.869604401\ldots$ to
$7.0666071493\ldots$.

We also consider the area ratio in \cref{conj:ratio}. By taking account
of the number of sides of each cell, we obtain the following improvement
of Fejes T\'oth's bound.

\begin{theorem}\label{thm:ratio}
There exists an absolute $n_0$ such that, for every $n\ge n_0$ and every
arrangement $\cA$ of $n$ distinct great circles on $S^2$, with no three
concurrent,
\[
 \frac{M(\cA)}{m(\cA)}>13.
\]
\end{theorem}

Both bounds are uniform over the arrangements, but
\cref{conj:ratio,conj:small-cell} remain open.

Related planar problems were studied by Barequet~\cite{Barequet01}
and Ismailescu~\cite{Ismailescu03}. In particular, Ismailescu proved
that $n$ lines whose pairwise intersections lie inside the unit disk
determine a cell of area at least $\pi/(4n)$. If no three lines are
concurrent, this gives a growing area ratio for the disk subdivision,
including boundary cells, but does not imply \cref{conj:ratio}.
For great circles, Bezdek and L\'angi~\cite{BL24} study extremal inradii
and circumradii. See also Maehara and Martini~\cite[\S4.7.3]{MM24}
for a discussion of Fejes T\'oth's tilings.

We begin in \cref{sec:kernel} with a perimeter inequality for planar
polygons and an estimate for the kernel arising from pairs of great
circles. In \cref{sec:transfer}, we use gnomonic projection to transfer
the perimeter bound to spherical cells. We prove \cref{thm:area} in
\cref{sec:smallest} by summing this bound over the cells. In
\cref{sec:ratio}, we also use a weighted estimate involving the side
counts to prove \cref{thm:ratio}. The numerical estimates are proved
in the appendix.

\section{Planar polygons and a spherical kernel}
\label{sec:kernel}

In this section, we prove the perimeter inequality which will be used
in both main results. We then show how the angles at a crossing give
a kernel of the normal vectors of two circles. The final lemma gives
a lower bound for the sum of a polynomial approximation to this kernel
over all pairs of circles.

For $0<t<\pi$, define
\[
 \psi(t)=\sqrt{\frac{2t\tan(t/2)}{\pi}}.
\]
We now state the planar inequality. Here an exterior angle is $\pi$
minus the interior angle. Recall that the exterior angles of a convex
planar polygon sum to $2\pi$, and that $p/\sqrt A$ is unchanged by
scaling the polygon.

\begin{lemma}\label{lem:polygon}
Let $P$ be a convex Euclidean $k$-gon with area $A$, perimeter $p$, and
exterior angles $t_1,\ldots,t_k\in(0,\pi)$. Then
\begin{equation}\label{eq:corner-planar}
 \frac{p}{\sqrt A}\ge\sum_{j=1}^k\psi(t_j).
\end{equation}
Moreover,
\begin{equation}\label{eq:polygonal-isoperimetric}
 \frac{p}{\sqrt A}\ge f(k),\qquad
 f(x)=2\sqrt{x\tan(\pi/x)}\quad(x>2).
\end{equation}
The function $f$ is decreasing and convex on $(2,\infty)$, and
\begin{equation}\label{eq:f-lower}
 f(x)\ge 2\sqrt\pi\left(1+\frac{\pi^2}{6x^2}\right).
\end{equation}
\end{lemma}

\begin{figure}[htbp]
\centering
\begin{tikzpicture}[font=\small,>=Stealth,line join=round]
\begin{scope}[xshift=0.0cm,scale=0.95]
\draw[thick,fill=gray!6] (1.15000,1.11881) -- (-0.04665,1.67682) -- (-1.51747,-0.42372) -- (-0.42364,-1.34155) -- (1.15000,-0.43301) -- cycle;
\draw[->,gray] (1.15000,0.34290)--(1.57000,0.34290);
\draw[->,gray] (0.55167,1.39782)--(0.72917,1.77847);
\draw[->,gray] (-0.78206,0.62655)--(-1.12610,0.86745);
\draw[->,gray] (-0.97055,-0.88264)--(-1.24052,-1.20438);
\draw[->,gray] (0.36318,-0.88728)--(0.57318,-1.25101);
\node at (0,-2.05) {$P$};\end{scope}
\begin{scope}[xshift=4.9cm,scale=0.95]
\draw[thick,fill=gray!6] (1.00000,0.63707) -- (-0.33786,1.26093) -- (-1.34474,-0.17704) -- (-0.10640,-1.21613) -- (1.00000,-0.57735) -- cycle;
\draw[dashed] (0,0) circle (1);
\fill (0,0) circle (1.2pt);
\draw (0,0)--node[above] {$1$} (1,0);
\draw[->,gray] (1.00000,0.02986)--(1.42000,0.02986);
\draw[->,gray] (0.33107,0.94900)--(0.50857,1.32965);
\draw[->,gray] (-0.84130,0.54194)--(-1.18534,0.78285);
\draw[->,gray] (-0.72557,-0.69658)--(-0.99554,-1.01832);
\draw[->,gray] (0.44680,-0.89674)--(0.65680,-1.26047);
\node at (0,-2.05) {$T$};\end{scope}

\node[align=center,font=\scriptsize] at (2.45,0.2) {same outward\\normals};
\draw[->] (1.8,-0.35)--(3.15,-0.35);
\begin{scope}[xshift=9.6cm,scale=1.25]
\coordinate (O) at (0,0);
\coordinate (V) at ({1/cos(35)},0);
\coordinate (U) at (35:1);
\coordinate (D) at (-35:1);
\draw[dashed] (0,0) circle (1);
\draw[thick] (U)--(V)--(D);
\draw (O)--(U) (O)--(D);
\draw[dotted] (O)--(V);
\draw (0.34,0) arc[start angle=0,end angle=35,radius=0.34];
\node at (18:0.60) {$t/2$};
\node[above left] at (U) {$U$};
\node[below left] at (D) {$D$};
\node[right] at (V) {$V$};
\node[left] at (O) {$O$};
\node[above left] at (35:0.55) {$1$};
\pic[draw,angle radius=4pt] {right angle=O--U--V};
\node[font=\scriptsize,align=center] at (0,-1.56)
 {$|VU|=|VD|=\tan(t/2)$};
\end{scope}
\end{tikzpicture}
\caption{The auxiliary polygon $T$ has the same outward normal
directions as $P$ and is circumscribed about the unit disk.
In the local construction on the right, $t$ is the exterior angle
at $V$. The right triangle $OUV$ has $|OU|=1$ and
$\angle UOV=t/2$, giving the tangent length used to compute $J$.}
\label{fig:tangent-polygon}
\end{figure}
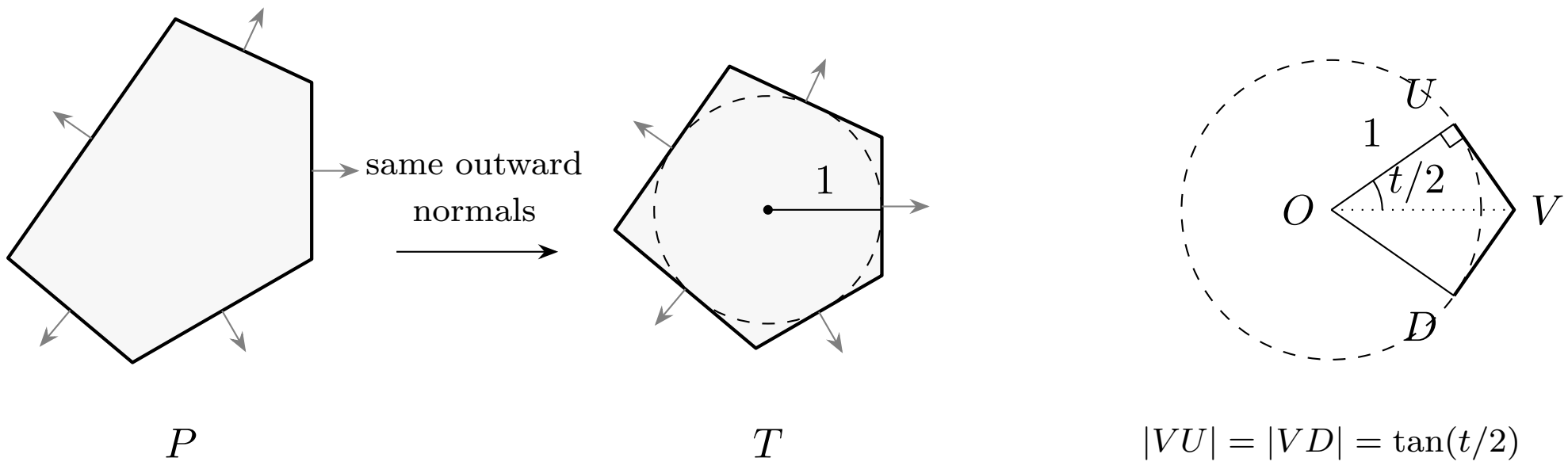

Note that both bounds are sharp for regular polygons. The first bound
depends on the individual angles, while the second depends only on
the number of sides. We will need both bounds in the proof of
\cref{thm:ratio}.

\begin{proof}
Let $T$ be the polygon circumscribed about the unit disk whose outward
edge normals agree with those of $P$; see \cref{fig:tangent-polygon}.
Thus we keep the side directions and move the supporting lines until
each is tangent to the unit disk. It is bounded because the normals
of a bounded polygon are not contained in a closed semicircle. Its area
is
\[
 J=\sum_{j=1}^k\tan(t_j/2).
\]
Indeed, the two tangent segments adjacent to a vertex each have length
$\tan(t_j/2)$, and the area is half the perimeter because the inradius
is one. Recall that the mixed area is defined by
\[
 \operatorname{area}(P+\varepsilon T)
 =A+2\varepsilon V(P,T)+\varepsilon^2J\qquad(\varepsilon\ge0),
\]
where $P+\varepsilon T=\{x+\varepsilon y:x\in P,\ y\in T\}$.
If $\ell_i$ and $\nu_i$ are the side lengths and outward unit normals
of $P$, the support-function formula gives
$2V(P,T)=\sum_i\ell_i h_T(\nu_i)=p$, since
$h_T(\nu_i)=\max_{x\in T}x\cdot\nu_i=1$ for every $i$.
The planar Minkowski inequality
$V(P,T)^2\ge A J$ therefore gives $p^2\ge4AJ$.
These mixed-area facts are standard consequences of the
Brunn--Minkowski inequality; see~\cite[Chapters~5 and~7]{Schneider14}.
Since $\sum_jt_j=2\pi$, Cauchy--Schwarz gives
\[
 \sum_{j=1}^k\psi(t_j)
 \le \sqrt{\frac2\pi
       \left(\sum_jt_j\right)\left(\sum_j\tan(t_j/2)\right)}
 =2\sqrt J,
\]
proving~\eqref{eq:corner-planar}. Jensen's inequality applied to
$\tan(t/2)$ gives $J\ge k\tan(\pi/k)$, which proves
\eqref{eq:polygonal-isoperimetric}.

To check monotonicity and convexity, put $z=\pi/x$, $T_0=\tan z$,
$S_0=\sec^2z$ and $h(x)=xT_0$. Then
\[
 h'=T_0-zS_0<0,\qquad h''=\frac{2z^2T_0S_0}{x}.
\]
The sign of $f''$ is the sign of
\[
 2hh''-(h')^2=4z^2T_0^2S_0-(zS_0-T_0)^2.
\]
Using $z\le T_0\le zS_0$, we obtain
$0\le zS_0-T_0\le zT_0^2$. The last expression is consequently at
least $z^2T_0^2(4S_0-T_0^2)>0$.
Finally, integrating $\tan'=1+\tan^2$ successively, starting with
$\tan z\ge z$, gives
\[
 \tan z\ge z+\frac{z^3}{3}+\frac{2z^5}{15}
          \ge z\left(1+\frac{z^2}{6}\right)^2.
\]
Taking square roots proves~\eqref{eq:f-lower}.
\end{proof}

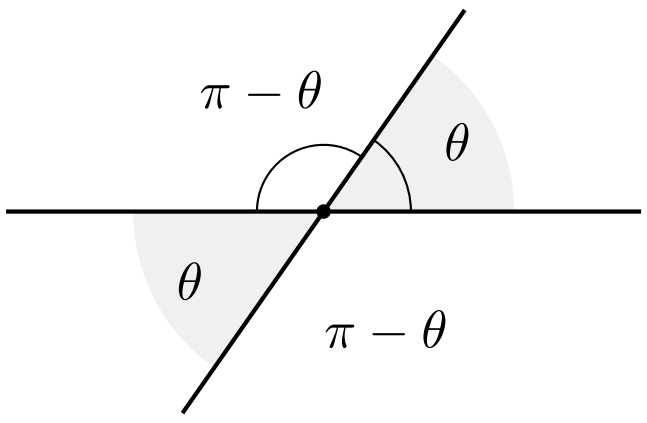
\begin{figure}[htbp]
\centering
\begin{tikzpicture}[scale=1.05,font=\small]
\fill[gray!12] (0,0)--(0:1.2) arc[start angle=0,end angle=55,radius=1.2]--cycle;
\fill[gray!12] (0,0)--(180:1.2) arc[start angle=180,end angle=235,radius=1.2]--cycle;
\draw[thick] (-2,0)--(2,0);
\draw[thick] (235:1.55)--(55:1.55);
\draw (0:0.55) arc[start angle=0,end angle=55,radius=0.55];
\draw (55:0.42) arc[start angle=55,end angle=180,radius=0.42];
\node at (27.5:0.95) {$\theta$};
\node at (117.5:0.85) {$\pi-\theta$};
\node at (207.5:0.95) {$\theta$};
\node at (297.5:0.85) {$\pi-\theta$};
\fill (0,0) circle (1.3pt);
\end{tikzpicture}
\caption{Tangent directions at one crossing of two great circles.
The labels show interior angles; the corresponding exterior angles
are their supplements, so the same two values each occur twice.
The antipodal crossing gives four copies of each value in total.}
\label{fig:crossing}
\end{figure}

At a crossing of two great circles, the exterior angles of the four
incident cells occur in complementary pairs; see \cref{fig:crossing}.
The next inequality is
the source of the constant in \cref{thm:area}.

\begin{lemma}\label{lem:kernel}
For $0<\theta<\pi$,
\begin{equation}\label{eq:kernel}
 \psi(\theta)+\psi(\pi-\theta)
 \ge b\,(1-\cos^2\theta)^{-1/4},
\end{equation}
where $b$ is defined in~\eqref{eq:constants}.
\end{lemma}

\begin{proof}
Put $x=\cos\theta$ and $y=(2/\pi)\arcsin x$. Direct substitution gives
\begin{align*}
 &(1-x^2)^{1/4}\bigl(\psi(\theta)+\psi(\pi-\theta)\bigr)\\
 &\quad=\sqrt{(1-x)(1-y)}+\sqrt{(1+x)(1+y)}\\
 &\quad=2\cos\!\left(\frac{\arcsin x-\arcsin y}{2}\right).
\end{align*}
The expression is even in $x$, so suppose $x\ge0$ and write
$a=\arcsin x\in[0,\pi/2]$. The function
\[
 d(a)=a-\arcsin(2a/\pi)
\]
is nonnegative because $\sin a\ge2a/\pi$ on $[0,\pi/2]$.
It vanishes at both endpoints and has derivative
\[
 d'(a)=1-\frac{2/\pi}{\sqrt{1-4a^2/\pi^2}}.
\]
Its unique maximum occurs at $a=a_*=\sqrt{\pi^2/4-1}$.
Thus $0\le d(a)\le d_*$, proving~\eqref{eq:kernel}.
\end{proof}

Before summing~\eqref{eq:kernel}, we need to deal with its singularity
when two normals become parallel. There are two reasons for this. The
sum over all ordered pairs of normals includes diagonal terms, where
the kernel is infinite. Also, the comparison with spherical polygons
will change the angles slightly, and we need to control the effect
of these changes. We therefore replace the kernel by a polynomial
lying below it, and replace $\psi$ by a Lipschitz function.
For an integer $L\ge0$, define
\begin{equation}\label{eq:truncation}
 \begin{gathered}
 a_j=\frac{(1/4)_j}{j!},\qquad
 P_L(z)=\sum_{j=0}^L a_jz^j,\qquad
 I_L=\sum_{j=0}^L\frac{a_j}{2j+1},\\
 B_L=\max\left\{\frac2\pi,\frac{2bP_L(1)}\pi\right\},\qquad
 \phi_L(t)=\min\{\psi(t),B_Lt\}.
 \end{gathered}
\end{equation}
Here $(1/4)_j$ denotes the rising factorial, with $(1/4)_0=1$.
We extend $\phi_L$ to $[0,\pi]$ by continuity; in particular,
$\phi_L(0)=0$ and $\phi_L(\pi)=B_L\pi$.
Since $\phi_L\le\psi$, replacing $\psi$ by $\phi_L$ preserves the
planar perimeter bound. The choice of $B_L$ ensures that clipping a
large angle does not destroy the lower bound for its complementary pair.

\begin{lemma}\label{lem:clipped}
The function $\phi_L$ is Lipschitz on $[0,\pi]$. It has a Lipschitz
constant $D_L=O(B_L^3)$, with an absolute implied constant. Moreover,
\begin{equation}\label{eq:clipped-kernel}
 \phi_L(\theta)+\phi_L(\pi-\theta)
 \ge bP_L(\cos^2\theta)\qquad(0<\theta<\pi).
\end{equation}
\end{lemma}

\begin{proof}
The quotient $\psi(t)/t$ is increasing, because $\tan(t/2)/t$ is
increasing. Since $\psi(\pi/2)/(\pi/2)=2/\pi\le B_L$, clipping can
occur only for $t\ge\pi/2$. If neither of the complementary angles is
clipped, \cref{lem:kernel} and the positive binomial expansion
\[
 (1-z)^{-1/4}=\sum_{j=0}^\infty a_jz^j\qquad(0\le z<1)
\]
give~\eqref{eq:clipped-kernel}. If clipping occurs, the clipped term
alone is at least $B_L\pi/2\ge bP_L(1)$. This also proves the claim.

For the derivative estimate, $\psi'$ is bounded on $[0,\pi/2]$ after
continuous extension at zero. On the unclipped part of $[\pi/2,\pi)$,
write $T_0=\tan(t/2)$. The inequality $\psi(t)\le B_Lt$ gives
$1\le T_0\le\pi^2B_L^2/2$, while $\psi(t)\ge\sqrt{T_0}$. Hence
\[
 \psi'(t)=\frac{T_0+(t/2)(1+T_0^2)}{\pi\psi(t)}=O(B_L^3).
\]
On the clipped part the derivative is $B_L$. There is a unique
transition point, so these bounds give the asserted Lipschitz constant.
\end{proof}

To apply the preceding lemma to an arrangement, choose a unit normal
$u_i$ to each great circle. Observe that the squared cosine of the
crossing angle is $(u_i\cdot u_h)^2$, regardless of the signs chosen
for the normals. We will use the following lemma to bound the sum
over all pairs. The lemma holds for any choice of unit vectors.

\begin{lemma}\label{lem:moments}
For unit vectors $u_1,\ldots,u_n\in S^2$ and every integer $j\ge0$,
\begin{equation}\label{eq:moments}
 \sum_{i,h=1}^n(u_i\cdot u_h)^{2j}\ge\frac{n^2}{2j+1}.
\end{equation}
Consequently,
\begin{equation}\label{eq:pair-moments}
 \sum_{i<h}P_L((u_i\cdot u_h)^2)
 \ge\frac12\bigl(n^2 I_L-nP_L(1)\bigr).
\end{equation}
\end{lemma}

\begin{proof}
Let $\mu=n^{-1}\sum_i\delta_{u_i}$ and let $\sigma$ be the uniform
probability measure on $S^2$. The kernel $(u\cdot v)^{2j}$ is positive
semidefinite: it is the inner product of the tensors $u^{\otimes2j}$
and $v^{\otimes2j}$. In particular, for any finite signed measure $\nu$,
\[
 \iint (u\cdot v)^{2j}\,d\nu(u)\,d\nu(v)
 =\left\|\int u^{\otimes2j}\,d\nu(u)\right\|^2\ge0.
\]
For fixed $u$, the coordinate $u\cdot v$ of a uniformly distributed
point $v\in S^2$ is uniform on $[-1,1]$. Thus the integral against
$\sigma$ in either variable is
\[
 \frac12\int_{-1}^1 s^{2j}\,ds=\frac1{2j+1}.
\]
Expanding its nonnegative energy on the signed measure $\mu-\sigma$
gives
\[
 0\le\iint(u\cdot v)^{2j}\,d(\mu-\sigma)(u)\,d(\mu-\sigma)(v)
 =\iint(u\cdot v)^{2j}\,d\mu(u)\,d\mu(v)-\frac1{2j+1}.
\]
This proves~\eqref{eq:moments}, including $j=0$. Multiply by $a_j>0$,
sum over $0\le j\le L$, and remove the $n$ diagonal terms, each equal
to $P_L(1)$. Every remaining unordered pair occurs twice, so division
by two gives~\eqref{eq:pair-moments}. This subtraction is legitimate
because the polynomial kernel is finite on the diagonal.
\end{proof}

\section{A uniform comparison with spherical cells}
\label{sec:transfer}

In order to apply \cref{lem:polygon}, we need a comparison between
planar and spherical cells. The main point is that the error must be
uniform, including for cells with many sides or very small angles.
We will first deal with cells whose perimeter is already large enough
for the required inequality. In the remaining case, the cell lies in
a small spherical cap, where we can control the distortion under
gnomonic projection.

For a cell $C$, write $A_C$, $p_C$ and $k_C$ for its area, perimeter and
number of sides. We use the identities
\begin{equation}\label{eq:arrangement}
 \begin{gathered}
 F=n^2-n+2,\qquad \sum_C k_C=4n(n-1)=4F-8,\\
 \sum_C p_C=4\pi n,\qquad \sum_C A_C=4\pi.
 \end{gathered}
\end{equation}
Indeed, there are $n(n-1)$ vertices and $2n(n-1)$ edges. Euler's formula
gives $F$, and each edge is incident with two cells. The perimeter sum
counts each great circle twice.

Every closed cell lies in an open hemisphere. To see this, choose signed
unit normals $v_1,\ldots,v_n$ so that its cone is
$\{x:v_i\cdot x\ge0\text{ for all }i\}$. For any nonzero point $x$ of
this cone,
\[
 \left(\sum_i v_i\right)\cdot x>0.
\]
Otherwise all $v_i\cdot x$ would vanish. The normals span $\R^3$
because no three great circles are concurrent, so this is impossible.
The cells are therefore geodesically convex polygons with $k_C\ge3$.

Let $t_{C,1},\ldots,t_{C,k_C}\in(0,\pi)$ be the exterior angles of $C$.
For a fixed truncation $L$, define
\begin{equation}\label{eq:cell-scores}
 \begin{gathered}
 q_C=\sum_{j=1}^{k_C}\phi_L(t_{C,j}),\qquad
 Q_L=2\pi B_L,\\
 s_C=\max\{q_C,f(k_C)\},\qquad
 U_L=\max\{Q_L,f(3)\}.
 \end{gathered}
\end{equation}
Gauss--Bonnet gives $\sum_jt_{C,j}=2\pi-A_C$. Since
$0\le\phi_L(t)\le B_Lt$, we have $0\le q_C\le Q_L$ and
$0\le s_C\le U_L$.
The quantity $q_C$ comes from the angles of $C$, while $f(k_C)$ depends
only on its number of sides. We take their maximum in order to use both
bounds from \cref{lem:polygon}. The following lemma gives the error
in passing to the sphere. We state it with a parameter $\alpha\le A_C$,
since we will later use both $\alpha=m(\cA)$ and $\alpha=A_C$.

\begin{lemma}\label{lem:transfer}
Let $C$ be a cell and let $0<\alpha\le A_C$ satisfy
$(U_L+1)\sqrt\alpha\le1$. Then
\begin{equation}\label{eq:transfer}
 p_C\ge\sqrt\alpha\,s_C
 -\frac{(U_L+1)^2}{4}\alpha^{3/2}(U_L+D_L k_C).
\end{equation}
The same inequality holds with $q_C$ in place of $s_C$.
\end{lemma}

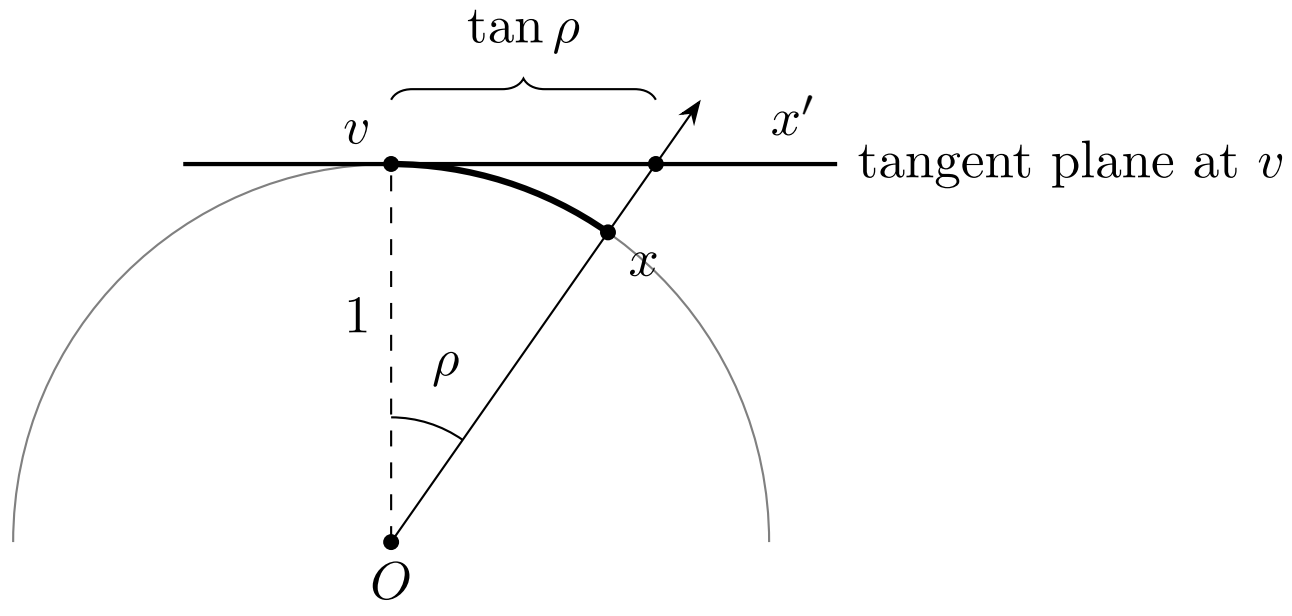
\begin{figure}[tbp]
\centering
\begin{tikzpicture}[scale=2.5,font=\small,>=Stealth]
\coordinate (O) at (0,0);
\coordinate (v) at (0,1);
\coordinate (x) at (55:1);
\coordinate (xp) at ({tan(35)},1);
\draw[gray] (-1,0) arc[start angle=180,end angle=0,radius=1];
\draw[thick] (-0.55,1)--(1.18,1);
\node[anchor=west] at (1.18,1) {tangent plane at $v$};
\draw[dashed] (O)--(v);
\draw[->] (O)--($(O)!1.17!(xp)$);
\draw[very thick] (55:1) arc[start angle=55,end angle=90,radius=1];
\draw (55:0.33) arc[start angle=55,end angle=90,radius=0.33];
\node at (72.5:0.48) {$\rho$};
\draw[decorate,decoration={brace,amplitude=4pt}]
 (0,1.17)--node[above=5pt] {$\tan\rho$} ({tan(35)},1.17);
\fill (O) circle (0.6pt) (v) circle (0.6pt)
 (x) circle (0.6pt) (xp) circle (0.6pt);
\node[below] at (O) {$O$};
\node[above left] at (v) {$v$};
\node[below right] at (x) {$x$};
\node[anchor=west] at (0.95,1.13) {$x'$};
\node[left] at (0,0.6) {$1$};
\end{tikzpicture}
\caption{A meridional cross-section of gnomonic projection on the unit
sphere. The ray from the center $O$ through $x$ meets the tangent
plane at $x'$. If the spherical distance from $v$ to $x$ is $\rho$,
then $|vx'|=\tan\rho$. The bold arc lies on the sphere.}
\label{fig:gnomonic}
\end{figure}

\begin{proof}
If $p_C\ge(U_L+1)\sqrt\alpha$, the assertion follows at once from
$s_C\le U_L$. Otherwise put $r=p_C/2<1/2$ and choose a vertex $v$ of
$C$. The shorter of the two boundary paths from $v$ to any other
vertex has length at most $r$. Thus all vertices lie in the spherical
cap of radius $r$ centered at $v$. The cell $C$ is the spherical convex
hull of its vertices. Since this cap is geodesically convex, it contains
the whole cell.

Gnomonic projection from the center of the sphere to the tangent plane
at $v$ maps $C$ to a convex Euclidean polygon $C'$ with the same number
of sides. Great-circle arcs become line segments because their planes
through the sphere's center intersect the tangent plane in lines.
At spherical distance $\rho$ from $v$, the projected radial distance
is $\tan\rho$; see \cref{fig:gnomonic}. Differentiating gives the radial
length factor $\sec^2\rho$, while the transverse factor is
$\tan\rho/\sin\rho=\sec\rho$ (with value one at $\rho=0$).
Consequently, if $A'$ and $p'$ are the area and perimeter of $C'$, then
\begin{equation}\label{eq:projection-metric}
 A'\ge A_C\ge\alpha,\qquad p'\le\sec^2r\,p_C.
\end{equation}

The condition number of the differential is at most $\sec r$. A linear
map with condition number $\kappa$ changes an angle by at most
$2\arcsin((\kappa-1)/(\kappa+1))$. For completeness, after orthogonal
changes of coordinates and a scalar dilation, such a map is
$\operatorname{diag}(1,\kappa)$. The maximum change in the direction
of one ray is $\arcsin((\kappa-1)/(\kappa+1))$, obtained by maximizing
$|\arctan(\kappa\tan u)-u|$; two rays give the stated bound.
It follows that each exterior angle changes by at most
\begin{equation}\label{eq:angle-distortion}
 2\arcsin\!\left(\tan^2(r/2)\right)\le r^2
 \qquad(0\le r\le1/2).
\end{equation}
The last inequality follows from $\tan(r/2)\le2r/3$ and
$\arcsin z\le9z/8$ for $0\le z\le1/9$.

Let $q'$ be the sum of $\phi_L$ over the exterior angles of $C'$, and
let $s'=\max\{q',f(k_C)\}$. By the Lipschitz bound and the fact that
the maximum is $1$-Lipschitz in its first argument,
\[
 q'\ge q_C-D_Lk_Cr^2,\qquad s'\ge s_C-D_Lk_Cr^2.
\]
By \cref{lem:polygon}, $p'\ge\sqrt{A'}s'$.
If $s_C-D_Lk_Cr^2\le0$, the right side
of~\eqref{eq:transfer} is nonpositive because
$r^2\le(U_L+1)^2\alpha/4$. Otherwise,
\begin{align*}
 p_C
 &\ge\cos^2r\,\sqrt\alpha\,(s_C-D_Lk_Cr^2)\\
 &\ge\sqrt\alpha\bigl(s_C-r^2(U_L+D_Lk_C)\bigr),
\end{align*}
using $\cos^2r\ge1-r^2$. Substituting
$r^2\le(U_L+1)^2\alpha/4$ proves~\eqref{eq:transfer}.
The bound with $q_C$ follows from $q_C\le s_C$.
\end{proof}

The local crossing count now converts the corner sum into a sum over
pairs of circles. Each pair contributes at two antipodal vertices,
which explains the factor four in the proof below.

\begin{lemma}\label{lem:total-score}
For every arrangement under consideration and every integer $L\ge0$,
\begin{equation}\label{eq:total-score}
 \sum_Cq_C\ge2b\bigl(n^2 I_L-nP_L(1)\bigr).
\end{equation}
\end{lemma}

\begin{proof}
Choose either unit normal $u_i$ for each circle and put
$\theta_{ih}=\arccos(u_i\cdot u_h)$. Distinct circles have
$0<\theta_{ih}<\pi$. At the two antipodal crossings of circles $i$
and $h$, the eight exterior corner angles are four copies each of
$\theta_{ih}$ and $\pi-\theta_{ih}$. Thus
\[
 \sum_Cq_C
 =4\sum_{i<h}\bigl(\phi_L(\theta_{ih})+
                         \phi_L(\pi-\theta_{ih})\bigr).
\]
Apply \cref{lem:clipped,lem:moments} to obtain~\eqref{eq:total-score}.
The estimate is independent of the chosen signs of the normals.
\end{proof}

\begin{remark}\label{rem:old-area}
We record the deduction of~\eqref{eq:old-area} from the classical
method. Fejes T\'oth~\cite[Lemma~3]{FT87} proves
\[
 p_C>\sqrt{A_C}\,f(k_C)-2A_C^{3/2}.
\]
Write $m=\min_CA_C$. Since $f$ is decreasing and convex,
\[
 \sum_C f(k_C)\ge F f(4-8/F)\ge4F.
\]
Also $\sum_C A_C^{3/2}\le(\sum_CA_C)^{3/2}=(4\pi)^{3/2}$.
Summing the perimeter inequality therefore gives
\[
 4\pi n\ge4F\sqrt m-2(4\pi)^{3/2}.
\]
Consequently $m\le(\pi^2+o(1))/n^2$, uniformly over the arrangements.
\end{remark}

\section{The smallest cell}
\label{sec:smallest}

We now prove \cref{thm:area}, which improves the upper bound for $H(n)$
in connection with \cref{conj:small-cell}. To see where the constant
comes from, temporarily ignore the error terms in the preceding
estimates. These give $4\pi n\ge\sqrt m\,2bI n^2$, and solving for $m$
gives the coefficient $C_*=(2\pi/(bI))^2$. We will make this argument
precise by choosing the truncation parameter as a function of $n$.
We first need the following estimates.

\begin{lemma}\label{lem:tails}
For $L\ge1$,
\begin{equation}\label{eq:tail-bounds}
 P_L(1)=O(L^{1/4}),\qquad
 I-I_L=O(L^{-3/4}),\qquad
 U_L=O(L^{1/4}),\qquad D_L=O(L^{3/4}),
\end{equation}
with absolute implied constants. Moreover, $I_L\uparrow I$.
\end{lemma}

\begin{proof}
The recurrence $a_j/a_{j-1}=1-3/(4j)$ gives
\[
 \log a_j\le-\frac34\sum_{h=1}^j\frac1h
             \le-\frac34\log(j+1).
\]
Hence $a_j\le(j+1)^{-3/4}$. The identity
$P_L(1)=(4L+1)a_L$, which also follows directly by induction, proves
the first bound. Integrating the positive binomial series and using
monotone convergence gives
\[
 I=\sum_{j=0}^\infty\frac{a_j}{2j+1},\qquad
 0<I-I_L=O\!\left(\sum_{j>L}j^{-7/4}\right)=O(L^{-3/4}).
\]
The remaining estimates follow from the definitions of $B_L$, $U_L$
and \cref{lem:clipped}.
\end{proof}

\begin{proof}[Proof of \cref{thm:area}]
Fix an arrangement of $n$ circles and put $m=m(\cA)$.
For any choice of $L$ satisfying $(U_L+1)\sqrt m\le1$, we may apply
\cref{lem:transfer} to each cell with $\alpha=m$. Summing and using
\eqref{eq:arrangement} and~\eqref{eq:total-score}, we obtain
\begin{align}
 4\pi n\ge\sqrt m\biggl[&2bn^2I_L-2bnP_L(1)\notag\\
 &-\frac{(U_L+1)^2}{4}
     m\bigl(FU_L+D_L(4F-8)\bigr)\biggr].
 \label{eq:master-area}
\end{align}
Since $mF\le4\pi$, the last term inside the brackets is at most
\[
 \pi(U_L+1)^2(U_L+4D_L).
\]
This bound is uniform over the arrangements, including cells with
arbitrarily many sides.

The three error bounds have orders $nL^{1/4}$, $L^{5/4}$ and
$n^2L^{-3/4}$, respectively. Choosing $L=n$ balances them at
$n^{5/4}$, below the leading order $n^2$.
By \cref{lem:tails} and the averaging bound for $m$,
\[
 (U_n+1)\sqrt m=O(n^{-3/4}),
\]
so the transfer condition holds for all sufficiently large $n$.
The diagonal term $nP_n(1)$, the curvature error
$(U_n+1)^2(U_n+4D_n)$, and the truncation error $n^2(I-I_n)$ are all
$O(n^{5/4})$. Thus~\eqref{eq:master-area} gives
\[
 4\pi n\ge\sqrt m\bigl(2bIn^2-O(n^{5/4})\bigr).
\]
The bracket is positive for large $n$. Dividing and squaring yields
\[
 m\le\frac{1}{n^2}
       \left(\frac{2\pi}{bI}+O(n^{-3/4})\right)^2
   =\frac{C_*+O(n^{-3/4})}{n^2}.
\]
Every implied constant is absolute. The estimate $C_*<7.067$ is
established in \cref{lem:constant-check} below.
\end{proof}

\section{The ratio of the largest and smallest cells}
\label{sec:ratio}

In this section, we prove \cref{thm:ratio}, giving an improved constant
bound for the ratio in \cref{conj:ratio}. We will use both the total
corner contribution and the polygonal bound $f(k_C)$. Suppose, for a
contradiction, that $M/m\le13$. We assign a weight in $[0,1]$ to each
cell according to its area, and let $t$ be the mean of these weights.
The total area determines $m$ in terms of $t$. We will show that the
resulting lower bound for the total perimeter exceeds $4\pi n$ for
every possible value of $t$. To do this, we first prove the following
elementary inequality.

\begin{lemma}\label{lem:polynomial}
For every $0\le t\le1$,
\begin{equation}\label{eq:scalar-gap}
 G(t):=\frac43+\frac{13}{5}t+
                \frac{533}{125}\frac{t^3}{(1+3t)^2}
       >\sqrt{1+12t}.
\end{equation}
\end{lemma}

\begin{proof}
Both sides are positive. On multiplying their squared difference by
$(1+3t)^4$, we obtain
\begin{align}
 P(t)&=(1+3t)^4\bigl(G(t)^2-(1+12t)\bigr)\notag\\
 &=\frac79+\frac{64}{15}t-\frac{301}{25}t^2
   -\frac{36416}{375}t^3-\frac{17977}{625}t^4
   +\frac{346908}{625}t^5+\frac{11957764}{15625}t^6.
 \label{eq:polynomial}
\end{align}
We prove that $P(t)>0$ for every $t\ge0$.

The signs of the coefficients of $P'$, in ascending order, are
$+,-,-,-,+,+$. Descartes' rule of signs gives at most two positive
roots. Since
\[
 P'(0)>0,\qquad
 P'(1/4)=-\frac{38603981}{6000000}<0,\qquad
 \lim_{t\to\infty}P'(t)=\infty,
\]
there are exactly two, one in $(0,1/4)$ and one in $(1/4,\infty)$.
Thus $P$ first increases and then decreases on $[0,1/4]$, and its
minimum there is attained at an endpoint. The endpoint values are
\[
 P(0)=\frac79,\qquad P(1/4)=\frac{27523969}{144000000},
\]
both greater than $19/100$.

For $t=1/4+u$ with $u\ge0$, the coefficients of $P(1/4+u)$, in
ascending powers of $u$, are
\[
 \frac{27523969}{144000000},\quad
 -\frac{38603981}{6000000},\quad
 \frac{7182083}{200000},\quad
 \frac{17257153}{37500},\quad
 \frac{17281483}{12500},\quad
 \frac{26609346}{15625},\quad
 \frac{11957764}{15625}.
\]
Comparison of coefficients gives
\[
 P(1/4+u)\ge\frac{19}{100}-\frac{161}{25}u+35u^2+460u^3.
\]
Since
\[
 u^3-\frac{3u}{625}+\frac{2}{15625}
 =\left(u-\frac1{25}\right)^2\left(u+\frac2{25}\right)\ge0,
\]
we conclude that
\[
 P(1/4+u)\ge35u^2-\frac{529}{125}u+\frac{1639}{12500}
          \ge\frac{1746}{546875}>0.
\]
This proves the lemma. In particular, the inequality is strict
uniformly on the compact interval $[0,1]$.
\end{proof}

\begin{proof}[Proof of \cref{thm:ratio}]
By \cref{lem:constant-check} in the appendix,
\[
 2bI>\frac{8\sqrt\pi}{3}.
\]
Since $I_L\uparrow I$, choose and fix $L$ such that
$2bI_L>8\sqrt\pi/3$. All constants depending on $L$ remain fixed for
the rest of this proof. By \cref{lem:total-score}, for all sufficiently
large $n$ we have
\begin{equation}\label{eq:score-margin}
 \frac1F\sum_Cq_C\ge\frac{8\sqrt\pi}{3}.
\end{equation}
Indeed, the left side is at least $2bI_L+O_L(n^{-1})$.

Suppose that an arrangement satisfies $M(\cA)/m(\cA)\le13$. Put
$m=m(\cA)$ and define
\[
 a_C=\frac{A_C}{m}\in[1,13],\qquad
 w_C=\frac{a_C-1}{12}\in[0,1],\qquad
 t=\frac1F\sum_Cw_C\in[0,1].
\]
Thus $w_C=0$ for a smallest cell, and $w_C=1$ would correspond to
area $13m$. These weights keep track of the extra perimeter forced
by areas larger than $m$. The chord bound for the concave square-root function gives
$\sqrt{a_C}\ge1+(\sqrt{13}-1)w_C$. By the definition of $s_C$,
\begin{equation}\label{eq:weighted-score}
 \sqrt{a_C}\,s_C
 \ge q_C+(\sqrt{13}-1)w_C f(k_C).
\end{equation}
This estimate makes no assumption about the relation between the areas,
angles and side counts of the cells.

Since $k_C\ge3$, $0\le w_C\le1$, and $\sum_Ck_C=4F-8$,
\[
 \sum_Cw_C k_C
 =3Ft+\sum_Cw_C(k_C-3)
 \le F(1+3t)-8.
\]
For $t>0$, the numbers $w_C/(Ft)$ sum to one, and the preceding
inequality bounds their weighted mean side count by $3+1/t$.
Jensen's inequality and the fact that $f$ is decreasing therefore give
\begin{equation}\label{eq:weighted-jensen}
 \sum_Cw_C f(k_C)\ge Ft\,f\!\left(3+\frac1t\right).
\end{equation}
When $t=0$, both sides are interpreted as zero; this agrees with the
continuous extension of the right side.

Every cell has $A_C\le13m\le52\pi/F$. Thus for all sufficiently large
$n$, \cref{lem:transfer} applies with $\alpha=A_C$. Its total error is
bounded by
\[
 \frac{(U_L+1)^2}{4}\sum_C A_C^{3/2}(U_L+D_Lk_C)
 =O_L(F^{-1/2})=O_L(n^{-1}).
\]
Moreover, the total area gives the exact identity
\[
 m=\frac{4\pi}{F(1+12t)}.
\]
Summing~\eqref{eq:weighted-score} and using
\eqref{eq:score-margin}--\eqref{eq:weighted-jensen}, we obtain
\begin{equation}\label{eq:ratio-master}
 4\pi n\ge\sqrt{4\pi F}\,
 \frac{8\sqrt\pi/3+(\sqrt{13}-1)t f(3+1/t)}{\sqrt{1+12t}}
 -O_L(n^{-1}).
\end{equation}

By~\eqref{eq:f-lower}, $\sqrt{13}-1>13/5$ and
$\pi^2/6>41/25$, the fraction in~\eqref{eq:ratio-master} is at least
\[
 2\sqrt\pi\,\frac{G(t)}{\sqrt{1+12t}}.
\]
By \cref{lem:polynomial} and continuity, there is an absolute
$\eta>0$ such that $G(t)/\sqrt{1+12t}\ge1+\eta$ for all $t\in[0,1]$.
It follows that
\[
 4\pi n\ge4\pi\sqrt F\,(1+\eta)-O_L(n^{-1}).
\]
Since $\sqrt F/n\to1$, this is impossible for all sufficiently large
$n$. All bounds are independent of the arrangement, proving the theorem.
\end{proof}

\appendix
\section{Elementary estimates for the constants}
\label{app:constants}

We finish by proving the numerical estimates used in the proofs of
\cref{thm:area,thm:ratio}. All comparisons below use rational numbers,
so the proofs do not depend on the displayed decimal approximations
to $b$, $I$ and $C_*$.

\begin{lemma}\label{lem:constant-check}
The constants in~\eqref{eq:constants} satisfy
\[
 b>1.9727,\qquad I>1.19813,\qquad
 2bI>\frac{8\sqrt\pi}{3},\qquad C_*<7.067.
\]
\end{lemma}

\begin{proof}
We use the elementary bounds $3.1415<\pi<3.1416$.
They follow, for example, from Machin's identity
$\pi=16\arctan(1/5)-4\arctan(1/239)$: four terms of the alternating
series for $\arctan(1/5)$, with its next-term error bound, and the first
two terms for $\arctan(1/239)$ suffice.

Let
\[
 c(x)=1-\frac{x^2}{2}+\frac{x^4}{24}-\frac{x^6}{720}.
\]
For $0<x<1$, the alternating cosine series gives $\cos x>c(x)$.
Direct rational comparisons give
\[
 (1.2114)^2>\frac{(3.1416)^2}{4}-1,\qquad
 c(0.8806)>\frac2{3.1415}.
\]
Since $\arcsin(2a_*/\pi)=\arccos(2/\pi)$, these inequalities imply
\[
 a_*<1.2114,\qquad
 \arccos(2/\pi)>0.8806,\qquad d_*<0.3308.
\]
Consequently
\[
 b>2\cos(0.1654)>2c(0.1654)>1.9727.
\]

To estimate $I$, substitute $x=1-u^4$. The positive binomial series
now converges geometrically and gives
\begin{align*}
 I&=4\,2^{-1/4}\int_0^1u^2(1-u^4/2)^{-1/4}\,du\\
  &=4\,2^{-1/4}\sum_{j=0}^\infty\frac{a_j}{2^j(4j+3)}.
\end{align*}
Taking the first ten terms, we have the exact identity
\[
 \sum_{j=0}^9\frac{a_j}{2^j(4j+3)}
   =\frac{57451188433862627}{161285831670104064}.
\]
The inequality $(0.840896)^4<1/2$ shows that
$2^{-1/4}>0.840896$. Thus
\[
 I>4(0.840896)
       \frac{57451188433862627}{161285831670104064}>1.19813.
\]
Finally, the rational comparisons
\[
 \bigl(2(1.9727)(1.198)\bigr)^2
       >\frac{64}{9}(3.1416)
\]
and
\[
 \bigl(2(3.1416)\bigr)^2
       <7.067\bigl((1.9727)(1.19813)\bigr)^2
\]
give the last two assertions. All terminating decimals in this proof
denote exact rational numbers.
\end{proof}

\section*{AI involvement in this research}
    The research presented in this paper has benefitted from AI collaboration, namely ChatGPT 6 Astra. OpenAI Codex has been used for brainstorming ideas and improving our draft bounds. ChatGPT is also used in generating some parts of the text.

    \section*{Acknowledgements} The author thank Alexander Polyanskii for emphasizing the significance of the problem considered in this paper. The author also thank P\'eter \'Agoston for helpful discussions. 
\bibliography{references}{}
\bibliographystyle{plain}

\end{document}